\documentclass[12pt]{amsart}

\usepackage{amsmath,amssymb}
\usepackage{amsfonts}
\usepackage{amsthm}
\usepackage{latexsym}
\usepackage{graphicx}
\usepackage{epsfig}

\usepackage{color}

\newtheorem{theorem}{Theorem}[section]

\newtheorem{prop}{Proposition}[section]
\newtheorem{rema}[prop]{Remark}

\makeatletter \@addtoreset{equation}{section} \makeatother

\def\ddt{\frac{d}{dt}}

\begin{document}

\title{The Backward behavior of the Ricci and
Cross Curvature Flows on $\mbox{SL}(2, \mathbb{R})$  }

\author{Xiaodong Cao
$^*$}
\thanks{$^*$Research partially supported by NSF grant DMS 0904432}

\address{Department of Mathematics,
  Cornell University, Ithaca, NY 14853-4201}
\email{cao@math.cornell.edu, jmg16@cornelll.edu,
lsc@math.cornell.edu}

\author{John Guckenheimer$^{\natural}$}
\thanks{$^{\natural}$Research
partially supported by DOE grant DE-FG02-93ER25164}

\author{Laurent Saloff-Coste$^{\flat}$}
\thanks{$^{\flat}$Research
partially supported by NSF grant DMS 0603866}


\renewcommand{\subjclassname}{%
  \textup{2000} Mathematics Subject Classification}

\date{Dec. 10,  2009}

\maketitle

\markboth{Xiaodong Cao, John Guckenheimer, Laurent Saloff-Coste}
{The Backward behavior of the Ricci and Cross Curvature Flows on
$\mbox{SL}(2, \mathbb{R})$ }

\begin{abstract}  This paper is concerned with properties of maximal
solutions of the Ricci and cross curvature flows on locally
homogeneous three-manifolds of type $\mbox{SL}_2(\mathbb R)$. We
prove that, generically, a maximal solution originates at a
sub-Riemannian geometry of Heisenberg type. This solves a problem
left open in earlier work by two of the authors.
\end{abstract}

\section{Introduction}
\subsection{Homogeneous evolution equations}
On a closed $3$-dimensional Riemannian manifold $(M,g)$, let $Rc$
be the Ricci tensor and $R$ be the scalar curvature. The
celebrated Ricci flow \cite{H3} starting from a metric $g_0$ is
the solution of
\begin{equation}\label{RF}
\left\{\begin{array}{l}\frac{\partial}{\partial
t}g=-2Rc\\
g(0)=g_0.\end{array}\right.\end{equation}

Another tensor, the cross curvature tensor, call it  $h$, is used in
\cite{chowhamiltonxcf} to define  the cross curvature flow (XCF)
on $3$-manifolds with either negative sectional curvature or
positive sectional curvature. In the case of negative
sectional curvature, the flow (XCF)
starting from a metric $g_0$ is the
solution of
\begin{equation}\label{XCF}
\left\{\begin{array}{l}\frac{\partial}{\partial
t}g=-2h\\
g(0)=g_0.\end{array}\right.
\end{equation}

Assume that computations are done in an orthonormal frame where
the Ricci tensor is diagonal. Then the cross curvature tensor is
diagonal and if the principal sectional curvatures are
$k_1,k_2,k_3$ ($k_i=K_{jkjk}$, circularly) and the Ricci and cross
curvature tensors are given by
\begin{equation}\label{R=kk}
R_{ii}=k_j+k_l \;\;(\mbox{circularly in $i,j,l$}),
\end{equation}
and
\begin{equation} \label{C=kk}
h_{ii} = k_jk_l.
\end{equation}\\

A very special case arises when the $3$-manifold $(M,g)$ is
locally homogeneous. In this case, both flows reduce to ODE
systems. As a consequence, the cross curvature flow can be defined
even if the sectional curvatures do not have a definite sign, as
this is the case for most of the homogeneous $3$-manifolds.

We say a Riemannian manifold $(M,g)$ is \textrm{locally
homogenous} if, for every two points $x\in M$ and $y\in M$, there
exist neighborhoods $U$ of $x$ and $V$ of $y$, and an isometry
$\varphi$ from $(U,g|_U)$ to $(V,g|_V)$ with $\varphi(x)=y$. We
say that $(M,g)$ is \textrm{homogeneous} if the isometry group
acts transitively on $M$, i.e., $U=V=M$ for all $x$ and $y$. By a
result of Singer \cite{singer60} , the universal cover of a
locally homogeneous manifold is homogeneous.

For a closed $3$-dimensional Riemannian manifold $(M, g_0)$ that
is locally homogeneous, there are 9 possibilities for the
universal cover. They can be labelled by the minimal isometry
group that acts transitively:
\begin{itemize}
\item[(a)]$H(3)$ ($H(n)$ denotes the isometry group of hyperbolic
$n$-space); $\mbox{SO}(3)\times \mathbb R$; $H(2)\times \mathbb
R$; \item[(b)] $\mathbb R^3$; $\mbox{SU}(2)$; $\mbox{SL}(2,\mathbb
R)$; Heisenberg; $E(1,1)=\mbox{Sol}$ (the group of isometry of
plane with a flat Lorentz metric); $E(2)$ (the group of isometries
of the Euclidean plane). This is called the Bianchi case in
\cite{isenberg1}.
\end{itemize}
The crucial difference between cases (a) and (b) above is that, in
case (b), the universal cover of the corresponding closed
$3$-manifold is (essentially) the minimal transitive group of
isometries itself (with the caveat that both $\mbox{SL}(2,\mathbb
R)$ and $E(2)$ should be replaced by their universal cover)
whereas in case (a) this minimal group is of higher dimension.
Case (b) corresponds exactly to the classification of
$3$-dimensional simply connected unimodular Lie groups
(non-unimodular Lie groups cannot cover a closed manifold). Any
$(M, g_0)$ in case (b) is of form $M=G/H$, where $G=\widetilde{M}$
is the universal cover, $H$ is a co-compact discrete subgroup of
$G$, and the metric $g_0$ descends from a left-invariant metric
$\widetilde{g_0}$ on $G$.

 The cases of $H(3)$,
$\mbox{SO}(3)\times \mathbb R$, $H(2)\times \mathbb R$ and
$\mathbb R^3$ all lead to well-understood and essentially trivial
behaviors for both the Ricci and cross curvature flows.

The forward behavior of the Ricci flow on locally homogeneous
closed $3$-manifolds was first analyzed by Isenberg and Jackson
\cite{isenberg1}. The forward and backward behaviors of the cross
curvature flows are treated in \cite{cnsc1,csc} whereas the
backward behavior of the Ricci flow is studied in \cite{cscbrf}.
Related works include \cite{gli08,knopf2,lot07}. In
\cite{cscbrf,csc}, the following interesting asymptotic behavior
of these Ricci and cross curvature flows in the backward direction
was observed: Let $g_t$ be a maximal solution defined on
$(-T_b,T_f)$ and passing through a generic $g_0$ at $t=0$. Then
either $T_b=\infty$ and $g(t)=e^{\lambda t}g_0$, or $T_b<\infty$
and there is a positive function $r(t)$ such that $r^2(t)g(t)$
converges to a sub-Riemannian metric of Heisenberg type, see
\cite{mon02,csc}. More precisely, in \cite{cscbrf,csc}, this
result was proved for all locally homogeneous closed
$3$-manifolds, except those of type $\mbox{SL}(2,\mathbb R)$.
Indeed, the structure of the corresponding ODE systems turns out
to be somewhat more complicated in the $\mbox{SL}(2,\mathbb R)$
case.\\

The aim of this paper is to prove the result described above in
the case of locally homogeneous $3$-manifolds of type
$\mbox{SL}(2,\mathbb R)$, i.e., $\widetilde{\mbox{SL}(2,\mathbb
R)}$. This will finish the proof of the following statement
announced in \cite{csc,cscbrf}.

\begin{theorem} Let $(M, g_0)$ be a complete locally homogeneous
$3$-manifold (compact or not), corresponding to the case (b)
discussed above.  Let $g(t), t\in (-T_b,T_f)$, be the maximal
solution of either the normalized Ricci flow {\em(\ref{normRF})}
or the cross curvature flow {\em (\ref{XCF})} passing through
$g_0$ at $t=0$. Let $d(t)$ be the corresponding distance function
on $M$. Assume that $g_0$ is generic among all locally homogeneous
metrics on $M$. Then
\begin{itemize}

\item either $M$ is of type $\mathbb R\sp 3$, $T_b=\infty$ and
$g(t)=g_0$,

\item or $T_b<\infty$ and there exists a function
$r(t):(-T_b,0)\to (0,\infty)$ such that, as $t$ tends to $-T_b$,
the metric spaces $(M,r(t)d(t))$ converge uniformly to a
sub-Riemannian metric space $(M,d(T_b))$ whose tangent cone at any
$m\in M$ is the Heisenberg group $\mathbb H_3$ equipped with its
natural sub-Riemannian metric.
\end{itemize}
\end{theorem}

Each of the manifolds $M$ considered in the above Theorem is of
the type $M=G/H$, where $G$ is a simply connected unimodular
$3$-dimensional Lie group and $H$ is a discrete subgroup of $G$.
By a locally homogeneous metric on $M$, we mean a metric that
descends from an invariant metric on $G$. A generic locally
homogeneous metric on $M$ is a metric that descends from a generic
invariant metric on $G$. Note that homogeneous metrics on $G$ can
be smoothly parametrized by an open set $\Omega$ in a finite
dimensional vector space and generic can be taken to mean ``for an
open dense subset of $\Omega$".

By definition, the uniform convergence of metric spaces $(M, d_t)$
to $(M,d)$ means the uniform convergence over compact sets of
$(x,y) \rightarrow d_t (x,y)$ to $(x,y) \rightarrow d (x,y)$, see
 \cite{csc,cscbrf} for  more details.

\subsection{The Ricci and cross curvature flows on
homogeneous $3$-manifolds}

Assume that $\mathfrak g$ is a $3$-dimensional real Lie unimodular
algebra equipped with an oriented Euclidean structure. According
to J. Milnor \cite{milnor76}
 there exists a (positively oriented) orthonormal basis
$(e_1,e_2,e_3)$ and reals $\lambda_1,\lambda_2,\lambda_3$ such
that the bracket operation of the Lie algebra has the form
$$[e_i,e_j]=\lambda_k e_k  \;\;\;\mbox{(circularly in $i,j,k$)}.$$
Milnor shows that such a basis diagonalizes the Ricci tensor and
thus also the cross curvature tensor. If $f_i= a_ja_k e_i$ with
nonzero $a_i,a_j,a_k\in \mathbb R$, then $[f_i,f_j]= \lambda_k
a_k^2 f_k$ (circularly in $i,j,k$). Using a choice of orientation,
we may assume that at most one of the $\lambda_i$ is negative and
then, the Lie algebra structure is entirely determined by the
signs (in $\{-1,0,+1\}$) of $\lambda_1,\lambda_2,\lambda_3$. For
instance, $+,+,+$ corresponds to $\mbox{SU}(2)$ whereas $+,+,-$
corresponds to $\mbox{SL}(2,\mathbb R)$.

In each case, let $\epsilon=(\epsilon_1,\epsilon_2,\epsilon_3)\in
\{-1,0,+1\}^3$ be the corresponding choice of signs. Then, given
$\epsilon$ and an Euclidean  metric $g_0$ on the corresponding Lie
algebra, we can choose a basis $f_1,f_2,f_3$ (with $f_i$ collinear
to $e_i$ above) such that
\begin{equation}\label{MF}
[f_i,f_j]= 2\epsilon_k f_k \;\;\;\mbox{(circularly in $i,j,k$)}.
\end{equation}
We call $(f_i)_1^3$ a Milnor frame for $g_0$. The metric, the
Ricci tensor and the cross curvature tensor are diagonalized in
this basis and this property is obviously maintained throughout
either the Ricci flow or cross curvature flow. If we let
$(f^i)_1^3$ be the dual frame of $(f_i)_1^3$, the metric $g_0$ has
the form
\begin{equation}\label{g_0}
g_0= A_0 f^1\otimes f^1 +B_0f^2\otimes f^2+ C_0f^3\otimes
f^3.\end{equation} Assuming existence of the flow $g(t)$ starting
from $g_0$, under either the Ricci flow or the cross curvature
flow (positive or negative), the original frame $(f_i)_1^3$ stays
a Milnor frame for $g(t)$ along the flow and $g(t)$ has the form
\begin{equation}\label{gflow}
g(t)= A(t) f^1\otimes f^1 +B(t) f^2\otimes f^2+ C(t) f^3\otimes
f^3.\end{equation} It follows that these flows reduce to ODEs in
$(A,B,C)$. Given a flow, the explicit form of the ODE depends on
the underlying Lie algebra structure. With the help of the
curvature computations done by Milnor in \cite{milnor76}, one can
find the explicit form of the equations for each Lie algebra
structure. The Ricci flow case was treated in \cite{isenberg1}.
The computations of the ODEs corresponding to the cross curvature
flow are presented in \cite{cnsc1,csc}.

\subsection{Invariant metrics on $\mbox{SL}(2,\mathbb R)$}

This paper is devoted to study of the Ricci and cross curvature
flows on $3$-dimensional Riemannian manifolds that are covered by
$\widetilde{\mbox{SL}(2,\mathbb R)}$. Since it makes no
differences, we focus on $\mbox{SL}(2,\mathbb R)$. Given a
left-invariant metric $g_0$ on $\mbox{SL}(2,\mathbb R)$, we fix a
Milnor frame $\{f_i\}_1^3$ such that
$$[f_2,f_3]=-2f_1, \;\;[f_3,f_1]=2f_2,\;\;[f_1,f_2]=2f_3$$
and $$g_0=A_0f^1\otimes f^1+B_0f^2\otimes f^2+C_0f^3\otimes f^3.$$
The sectional curvatures are
\begin{align*}
K(f_2 \wedge f_3)&=\frac{1}{ABC}(-3A^2+B^2+C^2-2BC-2AC-2AB),\\
K(f_3 \wedge f_1)&=\frac{1}{ABC}(-3B^2+A^2+C^2+2BC+2AC-2AB),\\
K(f_1 \wedge f_2)&=\frac{1}{ABC}(-3C^2+A^2+B^2+2BC-2AC+2AB).
\end{align*}

Recall that the Lie algebra $\mbox{sl}(2,\mathbb R)$ of
$\mbox{SL}(2,\mathbb R)$ can be realized as the space of two by
two real matrices with trace $0$. A basis of this space is
$$W=\left(\begin{array}{cc}0&-1\\1&0\end{array}\right),\;\;
H=\left(\begin{array}{cc}1&0\\0&-1\end{array}\right),\;\;
V=\left(\begin{array}{cc}0&1\\1&0\end{array}\right).$$ These
satisfy $$[H,V]=-2W,\; [W,H]=2V,\; [V,W]=2H.$$ This means that
$(W,V,H)$ can be taken as a concrete representation of the above
Milnor basis $(f_1,f_2,f_3)$. In particular, $f_1$ corresponds to
rotation in $\mbox{SL}(2,\mathbb R)$. Note further that exchanging
$f_2,f_3$ and replacing $f_1$ by $-f_1$ produce another Milnor
basis. This explains the $B,C$ symmetry of the formulas above.

\subsection{Normalizations}
Let $g(t), t\in I$, be a maximal solution of
\begin{equation}\label{RX}
\left\{\begin{array}{l}\frac{\partial}{\partial
t}g=-2v\\
g(0)=g_0,\end{array}\right.
\end{equation}
where $v$ denotes either the Ricci tensor $Rc$ or the cross curvature
tensor $h$. By
renormalization of $g(t)$, we mean a family
$\widetilde{g}(\widetilde{t}), \widetilde{t}\in \widetilde{I},$
obtained by a change of scale in space and a change of time, that
is
$$\widetilde{g}(\widetilde{t})=\psi (t) g(t),~
\widetilde{t}= \phi(t)$$
where $\phi$ is chosen appropriately. The choices of $\phi$ are
different for the two flows because of their different structures.
For the Ricci flow, take
$$\phi (t)=\int_{0}^t \psi(s) ds.$$
In the case of the cross curvature flow, take
$$\phi (t)=\int_{0}^t \psi^2(s) ds.$$
Now, set $\widetilde{\psi} (\widetilde{t})=\psi (t)$. Then we have
$$\frac{\partial \widetilde{g}}{\partial
\widetilde{t}}=-2\widetilde{v}+\left(\frac{d}{d \widetilde{t}} \ln
\widetilde{\psi}\right) \widetilde{g},$$ where $\widetilde{v}$ is
either the Ricci or the cross curvature tensor of $\widetilde{g}$.

On compact $3$-manifolds, it is customary to take $\ddt \ln
\psi=\frac{2}{3} \overline{v}$, where $\overline{v}=\frac{\int
tr(v) d\mu}{\int d\mu}$ is the average of the trace of either the
Ricci or the cross curvature tensor. In both cases, this choice
implies that the volume of the metric $\widetilde{g}$ is constant.
Obviously, studying any of the normalized versions is equivalent
to studying the original flow. Notice that the finiteness of $T_b$
or $T_f$ is not preserved under different normalization of flows.

\section{The Ricci Flow on $\mbox{SL}(2, \mathbb{R})$}

\subsection{The ODE system}

Mostly for historical reasons, we will consider the normalized
Ricci flow
\begin{equation}\label{normRF} \frac{\partial g}{\partial t}= - 2
Rc +\frac{2}{3}R g ,\;\; g(0)=g_0,\end{equation} where $g_0$ is a
left-invariant metric on $\mbox{SL}(2,\mathbb R)$. Let $g(t)$,
$t\in (-T_b,T_f)$ be the maximal solution of the normalized Ricci
flow through $g_0$.  In a Milnor frame $\{f_i\}_1^3 $ for $g_0$,
we write (see (\ref{gflow}))
$$g=  Af^1\otimes f^1+Bf^2\otimes f^2+Cf^3\otimes f^3.$$
Under (\ref{normRF}), $ABC=A_0B_0C_0$ is constant, and we set
$A_0B_0C_0 \equiv 4$. For this normalized Ricci flow, $A,B,C$
satisfy the equations
\begin{equation}\label{rpdesl2}
\left \{
\begin{aligned}
\frac{dA}{dt}=&\frac{2}{3}[-A^2(2A+B+C)+A(B-C)^2],\\
\frac{dB}{dt}=&\frac{2}{3}[-B^2(2B+A-C)+B(A+C)^2],\\
\frac{dC}{dt}=&\frac{2}{3}[-C^2(2C+A-B)+C(A+B)^2].
\end{aligned}
\right .
\end{equation}

\subsection{Asymptotic results}

Because of natural symmetries, we can assume without loss of
generality that $B_0 \geq C_0$. Then $B \geq C$ as long as a
solution exists. Throughout this section, we assume that $B_0\geq
C_0$.

\begin{theorem}[Ricci flow, forward direction, \cite{isenberg1}]
\label{th-fwsl2}
The forward time $T_f $ satisfies $T_f=\infty$. As $t$ tends to
$\infty$, $B-C$ tends to $0$ exponentially fast and
$$B(t)\sim (2/3)t , \;\; C(t)\sim (2/3)t \mbox{ and }A(t)\sim 9 t^{-2}.$$
\end{theorem}
In the backward direction, the following was proved in \cite{cscbrf}.
\begin{theorem}[Ricci flow, backward direction] \label{th-rfsl2}
We have $T_b\in (0,\infty)$, i.e., the maximal backward existence
time is finite. Moreover,
\begin{enumerate}
\item If there is a time $t<0$ such that $A(t)\ge B(t)$ then,
as $t$ tends to $-T_b$,
$$A(t)\sim \eta_1(t+T_b)^{-1/2},\;B(t)\sim \eta_2(t+T_b)^{1/4},\;
C(t)\sim \eta_3(t+T_b)^{1/4}$$ with $\eta_1=\sqrt{6}/4$  and
constants $\eta_i\in (0,\infty)$, $i=2,3$. \item If there is a
time $t<0$ such that $A\le B-C$ then, as $t$ tends to $-T_b$,
$$A(t)\sim \eta_1(t+T_b)^{1/4},\;B(t)\sim \eta_2(t+T_b)^{-1/2},\;
C(t)\sim \eta_3(t+T_b)^{1/4}$$
with $\eta_2=\sqrt{6}/4$, and constants $\eta_i\in (0,\infty)$, $i=1,3$.
\item If for all  time $t<0$, $B-C< A< B$ then,
as $t$ tends to $-T_b$,
$$A(t)\sim B(t)\sim \frac{\sqrt{6}}{4}(t+T_b)^{-1/2},\;
C(t)\sim \frac{32}{3}(t+T_b).$$
\end{enumerate}
\end{theorem}

As far as the normalized Ricci flow is concerned, the goal of his
paper is to show that the third case in the theorem above can only
occur when the initial condition $(A_0,B_0,C_0)$ belongs to a two
dimensional hypersurface. In particular, it does not occur for a
generic initial metric $g_0$ on $\mbox{SL}(2, \mathbb R)$.

\begin{theorem} \label{thm-2.3} Let $Q=\{(a,b,c)\in \mathbb R^3: a> 0, b>c>0\}$.
There is an open dense subset $Q_0$ of $Q$ such that,
for any maximal solution $g(t)$, $t\in (-T_b,T_f)$,
of the normalized Ricci flow with initial condition
$(A(0),B(0),C(0)\in Q_0$, as $t$ tends to $-T_b$,
\begin{enumerate}
\item either
$A(t)\sim (\sqrt{6}/4)(t+T_b)^{-1/2},\;B(t)\sim \eta_2(t+T_b)^{1/4},\;
C(t)\sim \eta_3(t+T_b)^{1/4}$
\item or
$A(t)\sim \eta_1(t+T_b)^{1/4},\;B(t)\sim (\sqrt{6}/4)(t+T_b)^{-1/2},\;
C(t)\sim \eta_3(t+T_b)^{1/4}$
\end{enumerate}
In fact, let $Q_1$ (resp. $Q_2$) be the set of initial conditions
such that case {\em (1)} (resp. case {\em (2)}) occurs. Then there
exists a smooth embedded hypersurface $S_0\subset Q$ such that
$Q_1,Q_2$ are the two connected components of $Q\setminus S_0$.
Moreover, for initial condition on $S_0$, the behavior is given by
case {\em (3)} of {\em Theorem \ref{th-rfsl2}}.
\end{theorem}

In order to prove this result, it suffices to study case (3) of
Theorem \ref{th-rfsl2}. This is done in the next section by
reducing the system (\ref{rpdesl2}) to a $2$-dimensional system.

\begin{rema}\label{rmk2.1}
The study below shows that, when the initial condition varies, all
values larger than $1$ of the ratio $\eta_2/\eta_3$ are attained
in case (1). Similarly, as the initial condition varies, all
positive values of the ratio $\eta_3/\eta_1$ are attained in case
(2).
\end{rema}

\subsection{The two-dimensional ODE system for the
Ricci flow}

For convenience, we introduce the backward normalized Ricci flow,
for which the ODE is
\begin{equation}\label{-rpdesl2}
\left \{
\begin{aligned}
\frac{dA}{dt}=&-\frac{2}{3}[-A^2(2A+B+C)+A(B-C)^2],\\
\frac{dB}{dt}=&-\frac{2}{3}[-B^2(2B+A-C)+B(A+C)^2],\\
\frac{dC}{dt}=&-\frac{2}{3}[-C^2(2C+A-B)+C(A+B)^2].
\end{aligned}
\right .
\end{equation}
By Theorem \ref{th-rfsl2}, the maximal forward solution of this
system is defined on $[0, T_b)$ with $T_b<\infty$.

 We start with the obvious
observation that if $(A, B, C)$ is a solution, then
$$t\mapsto (\lambda A (\lambda^2 t), \lambda B (\lambda^2 t),
\lambda C (\lambda^2 t))$$ is also a solution. By Theorem
\ref{th-rfsl2}, there are solutions with initial values in $Q$
such that $B/A$ tends to $\infty$ and others such that $B/A$ tends
to $0$. Let $Q_1$ be the set of initial values in $Q$ such that
$B/A$ tends to $0$ and $Q_2$ be the set of those for which $B/A$
tends to $\infty$. Let $S_0$ be the complement of $Q_1 \cup Q_2$
in $Q$. Again, by Theorem \ref{th-rfsl2}, for initial solution in
$S_0$, $B/A$ tends to $1$.

The sets $Q_1$, $Q_2$, $S_0$ must be homogeneous cones, i.e., are
preserved under dilations. Hence, they are determined by their
 projectivization on the plane $A=1$. So we set
$b=B/A,\;\;c=C/A$ and compute

\begin{equation}\label{ode-Abc}
\left\{\begin{array}{ccc}d b/dt &=& 2A^{2}b(1+b)(b-c-1)\\
dc/dt &=& -2 A^{2} c(1+c)(b-c+1).\end{array}\right.
\end{equation} This means that, up to a monotone time change,
the stereographic projection of any
flow line of (\ref{-rpdesl2}) on the plane $A=1$ is a flow line of
the planar ODE system
\begin{equation}\label{ode-bc}
\left\{\begin{array}{ccc}d b/dt &=& b(1+b)(b-c-1)\\
dc/dt &=& -c(1+c)(b-c+1).\end{array}\right.
\end{equation}
Set $\Omega=\{b > c > 0 \}$. By Theorem \ref{th-rfsl2}, any
integral curve of (\ref{ode-bc}) tends in the forward time
direction to either $(0,0)$, $(\infty, c_{\infty})$,
$0<c_{\infty}<\infty$, or $(1,0)$. The equilibrium points of
(\ref{ode-bc}), i.e., the points where $( db/dt,dc/dt)=(0,0)$,
 are $(1,0),(0,0)$ (notice that they are in $\partial \Omega$ but not in
 $\Omega$).
To investigate the nature of these equilibrium points, we compute
the Jacobian of the right-hand side of the (\ref{ode-bc}) which is
$$\left(\begin{array}{cc} 3b^2-3bc-c-1& -b^2-b\\
-c^2-c & 3c^2-2bc-b-1 \end{array}\right).$$

In particular, at $(0,0)$, this is $-1$ times the identity matrix
and the equilibrium point $(0,0)$ is attractive.  At $(1,0)$, the
Jacobian is $\left(\begin{array}{cc}2&-2\\0&-2\end{array}\right)$.
This point is a hyperbolic saddle point (the eigenvalues are $2$
and $-2$). By Theorem \ref{th-rfsl2}, any integral curve of
(\ref{ode-bc}) ending at $(1,0)$ must stay in the region
$\{b-c<1<b\}$. In that region, $b$ and $c$ are decreasing
functions of time $t$ and $\frac{ d c}{d b}$ is positive. Using
this observation and the stable manifold theorem, we obtain a
smooth increasing function $\phi : [1, \infty) \rightarrow
[0,\infty)$ whose graph $\gamma=\{(b,c): c=\phi (b)\}$ is  the
stable manifold at $(1,0)$ in $\overline{\Omega}$. By Theorem
\ref{th-fwsl2} and (\ref{ode-bc}), $\gamma$ is asymptotic to $c=b$
at infinity, and $\phi' (1)=2$. In particular, $\Omega\setminus
\gamma$ has two components $\Omega_1$ and $\Omega_2$ where $(0,0)
\in \overline{\Omega_1}$.  Further any initial condition in
$\Omega$ whose integral curve tends to $(1,0)$ must be on
$\gamma$. It is now clear that the cases (1), (2) and (3) in
Theorem \ref{th-rfsl2} correspond respectively to initial
conditions in $\Omega_1$, $\Omega_2$ and $\gamma$. This proves
Theorem \ref{thm-2.3} with $Q_i$ the positive cone with base
$\Omega_i$ and $S_0$ the positive cone with base $\gamma$.

\begin{figure} [!hbp]
   \begin{center}
   \begin{tabular}{c}
   \epsfig{figure=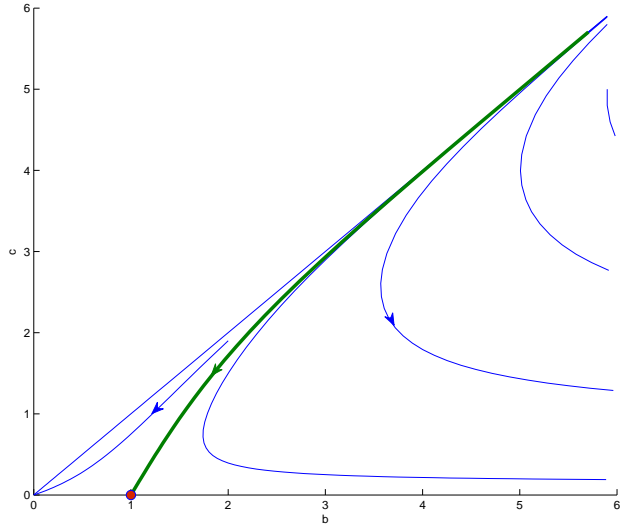,width=10cm}
   \end{tabular}
   \end{center}
   \caption
   { \label{fig:ricci}
   The flow line diagram of system (\ref{ode-bc})}
\end{figure}

Figure (\ref{fig:ricci}) shows the curve $\gamma$ and some flow
lines of (\ref{ode-bc}). It is easy to see from (\ref{ode-bc})
that the flow lines to the right of $\gamma$ have horizontal
asymptotes $\{c=c_{\infty}\}$ and that all positive values of
$c_{\infty}$ appear. This proves Remark \ref{rmk2.1} in case (2)
of Theorem \ref{thm-2.3}. The proof in case (1) is similar, but a
different choice of coordinates must be made.

\section{The cross curvature flow on $\mbox{SL}(2,\mathbb R)$}
\subsection{The ODE system}

We now consider the cross curvature flow (\ref{XCF}) where $g_0$
is a left-invariant metric on $\mbox{SL}(2,\mathbb R)$. Let
$g(t)$, $t\in (-T_b,T_f)$, be the maximal solution of the cross
curvature flow through $g_0$. Writing
$$g=  Af^1\otimes f^1+Bf^2\otimes f^2+Cf^3\otimes f^3,$$ we obtain
the system
\begin{equation}\label{pdesl2}
\left \{
\begin{aligned}
\frac{dA}{dt}=&-\frac{2AF_2F_3}{(ABC)^2},\\
\frac{dB}{dt}=&-\frac{2BF_3F_1}{(ABC)^2},\\
\frac{dC}{dt}=&-\frac{2CF_1F_2}{(ABC)^2},
\end{aligned}
\right .
\end{equation}
where
\begin{align*}
F_1=&-3A^2+B^2+C^2-2BC-2AC-2AB, \\
F_2=&-3B^2+A^2+C^2+2BC+2AC-2AB, \\
F_3=&-3C^2+A^2+B^2+2BC-2AC+2AB.
\end{align*}

\subsection{Asymptotic Results}

Without loss of generality we assume throughout this section that
$B_0 \geq C_0$. Then $B \geq C$ as long as a solution exists. The
behavior of the flow in the forward direction can be summarized as
follows. See \cite{cnsc1}.
\begin{theorem}[(XCF) Forward direction]\label{thm:fxcf}
If $B_0=C_0$ then $T_f=\infty$ and there exists a constant
$A_\infty \in (0,\infty)$ such that,
$$B(t)=C(t)\sim (24 A_\infty t)^{1/3} \mbox{ and }A(t)
\sim A_\infty \; \mbox{ as } t\rightarrow \infty.$$ If $B_0>C_0$,
$T_f$ is finite and there exists a constant $E\in (0,\infty)$ such
that,
$$ A(t)\sim B(t)\sim E(T_f-t)^{-1/2} \mbox{ and } C(t)\sim
8(T_f-t)^{1/2}\; \mbox{ as } t\rightarrow T_f.$$
\end{theorem}

In the backward direction, the following was proved in \cite{csc}.
\begin{theorem}[(XCF) Backward direction] \label{th-xcfsl2}
We have $T_b\in (0,\infty)$, i.e., the maximal backward existence
time is finite. Moreover,
\begin{enumerate}
\item If there is a time $t<0$ such that $A(t)\ge B(t)-C(t)$ then,
as $t$ tends to $-T_b$,
$$A(t)\sim \eta_1(t+T_b)^{-1/14},\;B(t)\sim \eta_2(t+T_b)^{3/14},\;
C(t)\sim \eta_3(t+T_b)^{3/14}$$
for some constants $\eta_i\in (0,\infty)$, $i=1,2,3$.
\item If there is a time $t<0$ such that
$A<\frac{1}{3}\left(2\sqrt{B^2-BC+C^2} -B-C\right)$ then,
as $t$ tends to $-T_b$,
$$A(t)\sim \eta_1(t+T_b)^{3/14},\;B(t)\sim \eta_2(t+T_b)^{-1/14},\;
C(t)\sim \eta_3(t+T_b)^{3/14}$$
for some constants $\eta_i\in (0,\infty)$, $i=1,2,3$.
\item If for all  time $t<0$,
$\frac{1}{3}\left(2\sqrt{B^2-BC+C^2} -B-C\right)\le A< B-C$ then,
as $t$ tends to $-T_b$,
$$A(t)\sim 64\eta^{-1}(t+T_b),\;B(t)\sim \eta +4\sqrt{2(t+T_b)},\;
C(t)\sim \eta -4\sqrt{2(t+T_b)}$$
for some  $\eta\in (0,\infty)$.
\end{enumerate}
\end{theorem}

As far as the cross curvature flow is concerned, the goal of this
paper is to show that the third case in the theorem above can only
occur when the initial condition $(A_0,B_0,C_0)$ belongs to a two
dimensional hypersurface. In particular, it does not occur for a
generic initial metric $g_0$ on $\mbox{SL}(2,\mathbb R)$.
\begin{theorem} \label{thm3.3} Let $Q=\{(a,b,c)\in \mathbb R^3: a> 0, b>c>0\}$.
There is an open dense subset $Q_0$ of $Q$ such that, for any
maximal solution $g(t)$, $t\in (-T_b,T_f)$, of the cross curvature
flow with initial condition $(A(0),B(0),C(0))\in Q_0$, as $t$
tends to $-T_b$,
\begin{enumerate}
\item either
$A(t)\sim \eta_1(t+T_b)^{-1/14},\;B(t)\sim \eta_2(t+T_b)^{3/14},\;
C(t)\sim \eta_3(t+T_b)^{3/14}$
\item or
$A(t)\sim \eta_1(t+T_b)^{3/14},\;B(t)\sim \eta_2(t+T_b)^{-1/14},\;
C(t)\sim \eta_3(t+T_b)^{3/14}$.
\end{enumerate}
 In fact, let $Q_1$ (resp. $Q_2$) be the set
of initial conditions such case {\em (1)} (resp. case {\em (2)}) occurs.
Then there exists a smooth embedded
hypersurface $S_0\subset Q$ such that $Q_1,Q_2$ are the two connected
components of $Q\setminus S_0$. Moreover,
for initial condition on $S_0$, the behavior is given by case
{\em (3)} of {\em Theorem \ref{th-xcfsl2}}.
\end{theorem}

In order to prove this result, it suffices to study case (3)
of Theorem \ref{th-xcfsl2}. In that case, it is proved in \cite{csc}
that $A$, $B$ and $C$ are monotone ($A,B$ non-decreasing,
$C$ non-increasing) on $(-T_b,0]$. In order to understand the behavior
of the solution, and because of the homogeneous structure of the
ODE system (\ref{pdesl2}), we can pass to the affine coordinates
$(A/C,B/C)$. This leads to a two-dimensional ODE system whose orbit
structure can be analyzed.

\begin{rema}\label{rmk3.1}
The analysis below shows that, in case (1) of Theorem \ref{thm3.3}
and when the initial condition varies, all the values larger than
$1$ of the ratio $\eta_2/\eta_3$ are attained. Similarly, in case
(2), all the values of the ratio $\eta_1/\eta_3$ are attained.
\end{rema}

\subsection{The two-dimensional ODE system for the cross
curvature flow}

For convenience, we introduce the backward cross curvature flow,
for which the ODE is
\begin{equation}\label{-xcpdesl2}
\left \{
\begin{aligned}
\frac{dA}{dt}=&\frac{2AF_2F_3}{(ABC)^2},\\
\frac{dB}{dt}=&\frac{2BF_3F_1}{(ABC)^2},\\
\frac{dC}{dt}=&\frac{2CF_1F_2}{(ABC)^2},
\end{aligned}
\right .
\end{equation} where $\{F_i\}_1^3$ are defined as before.
By Theorem \ref{th-xcfsl2}, the maximal forward solution of this
system is defined on $[0, T_b)$ with $T_b<\infty$.

Note that if $(A, B, C)$ is a solution, then
$$t\mapsto (\lambda A (t/\lambda^2), \lambda B (t/\lambda^2),
\lambda C (t/\lambda^2))$$ is also a solution. By Theorem
\ref{th-xcfsl2}, there are solutions with initial values in $Q$
such that $A/B$ tends to $\infty$ and others such that $(A/B,
C/B)$ tends to $(0, 0)$. Let $Q_1$ be the set of initial values in
$Q$ such that $A/B$ tends to $\infty$ and $Q_2$ be the set of
those for which $(A/B, C/B)$ tends to $(0, 0)$. Let $S_0$ be the
complement of $Q_1 \cup Q_2$ in $Q$. Again, by Theorem
\ref{th-xcfsl2}, for initial solution in $S_0$, $(A/B, C/B)$ tends
to $(0, 1)$.

The sets $Q_1$, $Q_2$, $S_0$ must be homogeneous cones, i.e., are
preserved under dilations. Hence, they are determined by their
projectivization on the plane $B=1$. So we set $a=A/B,\;\;c=C/B$
and compute

\begin{equation}\label{ode-Bac}
\left\{\begin{array}{ccc}d a/dt &=& \frac{8}{(Bac)^{2}}a(a+1)(a+c-1)
 \phi_3\\
dc/dt &=& -\frac{8}{(Bac)^2} c(1-c)(a+c+1)
\phi_1,\end{array}\right.
\end{equation} where
\begin{align*}
\phi_1=&-3a^2+1+c^2-2c-2ac-2a, \\
\phi_3=&-3c^2+a^2+1+2c-2ac+2a.
\end{align*}

This means that, up to a monotone time change, the stereographic
projection of any flow line of (\ref{-xcpdesl2}) on the plane
$B=1$ is a flow line of the planar ODE system
\begin{equation}\label{ode-ac}
\left\{\begin{array}{ccc}d a/dt &=& a(a+1)(a+c-1)
 \phi_3\\
dc/dt &=& -c(1-c)(a+c+1) \phi_1.\end{array}\right.
\end{equation}
Set $\Omega=\{a>0, 1 > c > 0 \}$. By Theorem \ref{th-xcfsl2}, any
integral curve of (\ref{ode-Bac}) tends in the forward time
direction to either $(0,0)$, $(\infty, c_{\infty})$,
$0<c_{\infty}<\infty$, or $(0,1)$. The equilibrium points of
(\ref{ode-Bac}), i.e., the points where $( da/dt,dc/dt)=(0,0)$,
 are $(0,0),(1,0), (0,1)$.
To investigate the nature of these equilibrium points, we compute
the Jacobian of the right-hand side of the ODE which is
$$\left(\begin{array}{cc} Y_{11}&
 Y_{12}\\
Y_{21} & Y_{22}
\end{array}\right),$$
where
\begin{eqnarray*}
Y_{11}&=&(3a^2+2ac+c-1)\phi_3+2a(a+1)(a+c-1)(a-c+1),\\
Y_{12}&=&a(a+1)[\phi_3+2(a+c-1)(-3c+1-a)],\\
Y_{21}&=&-c(1-c)[\phi_1-2(a+c+1)(3a+c+1)],\\
Y_{22}&=&(3c^2+2ac-a-1)\phi_1-2c(1-c)(a+c+1)(c-1-a).
\end{eqnarray*}

In particular, at $(0,0)$, this is $-1$ times the identity matrix
and the equilibrium point $(0,0)$ is attractive. Its basin of
attraction corresponds to region $Q_2$ defined above.  At $(1,0)$,
the Jacobian is
$\left(\begin{array}{cc}8&8\\0&8\end{array}\right)$. This point is
a repelling fixed point. It reflects the behavior of the forward
Ricci flow described in Theorem \ref{thm:fxcf}($B>C$).  At
$(0,1)$, the Jacobian is
$\left(\begin{array}{cc}0&0\\0&0\end{array}\right)$. This is the
equilibrium point of interest to us and a more detailed analysis
is required to determine trajectories that tend toward it. This is
done with coordinate transformations that \emph{blow-up} the
equilibrium.

Blow-up transformations introduce coordinates in which the blown
up equilibrium becomes a circle or projective line representing
directions through the equilibrium. Blow-up transformations reduce
the analysis of flows near degenerate equilibria to flows with
less degenerate equilibria. The blown up system allows us to
analyze the trajectories that are asymptotic to the equilibrium.
In particular, trajectories approaching the equilibrium of the
original system from different directions yield different
equilibria in the blown up system.

Translating the equilibrium point to the origin by setting
$e=c-1$, the equations (\ref{ode-Bac}) become
\begin{equation}\label{ode-Bae}
\left\{\begin{array}{ccc}
da/dt &=& -a(a+1)(a+e)(3e^2+4e+2ae-a^2), \\
de/dt &=& e(e+1)(a+e+2)(e^2-2ae-4a-3a^2).
\end{array}\right.
\end{equation}
The leading order terms of equations (\ref{ode-Bae}) have degrees
3 and 2: the next coordinate transformation $a=u^2, e = v$ of the
region $a \ge 0$ produces a system in which the leading terms of
both equations have degree 3:
\begin{equation}\label{ode-uv}
\left\{\begin{array}{ccc}
du/dt &=& -\frac{1}{2}u(u^2+1)(v+u^2)(3v^2+4v+2u^2v-u^4), \\
dv/dt &=& v(v+1)(v+2+u^2)(v^2-2u^2v-4u^2-3u^4).
\end{array}\right.
\end{equation}
To blow up the origin, the equations (\ref{ode-uv}) are
transformed to polar coordinates $(u,v) =
(r\cos(\theta),r\sin(\theta))$ and then rescaled by a common
factor of $r$, yielding the vector field $X$ defined by
\begin{equation}\label{ode-rth}
\left\{\begin{array}{ccl} dr/dt &=& 1/2\,r [ -33\,{r}^{3} \left(
\cos \left( \theta \right)  \right) ^{4} \left( \sin \left( \theta
\right)  \right) ^{3}-29\,{r}^{2}
 \left( \cos \left( \theta \right)  \right) ^{4} \left( \sin \left(
\theta \right)  \right) ^{2} \\
& & -11\,{r}^{4} \left( \cos \left( \theta
 \right)  \right) ^{6} \left( \sin \left( \theta \right)  \right) ^{2}
-{r}^{5} \left( \cos \left( \theta \right)  \right) ^{8}\sin
\left( \theta \right) -5\,{r}^{3} \left( \cos \left( \theta
\right)  \right)
^{6}\sin \left( \theta \right) \\
& & + {r}^{6} \left( \cos \left( \theta
 \right)  \right) ^{10} -35\,r \left( \cos \left( \theta \right)
 \right) ^{2} \left( \sin \left( \theta \right)  \right) ^{3}
 -20\, \left( \cos \left( \theta \right)  \right) ^{2} \left( \sin \left(
\theta \right)  \right) ^{2} \\
& &-4\,r \left( \cos \left( \theta \right)
 \right) ^{4}\sin \left( \theta \right) +{r}^{4} \left( \cos \left(
\theta \right)  \right) ^{8} +2\,{r}^{2} \left( \sin \left( \theta
 \right)  \right) ^{6}\\
& & -2\,{r}^{3} \left( \sin \left( \theta \right)
 \right) ^{5} \left( \cos \left( \theta \right)  \right) ^{2}-18\,{r}^
{2} \left( \cos \left( \theta \right)  \right) ^{2} \left( \sin
 \left( \theta \right)  \right) ^{4} \\
& & -10\,{r}^{4} \left( \sin \left( \theta \right)  \right) ^{4}
\left( \cos \left( \theta \right)
 \right) ^{4}
 + 6\,r \left( \sin \left( \theta \right)  \right) ^{5} \\
& & -6\,{r}^{5} \left( \sin \left( \theta \right)  \right) ^{3}
\left( \cos
 \left( \theta \right)  \right) ^{6}+4\, \left( \sin \left( \theta
 \right)  \right) ^{4}],
\\
d\theta/dt &=& -1/2\,\cos \left( \theta \right) \sin \left( \theta
\right)  [ -2 \,{r}^{2} \left( \sin \left( \theta \right)  \right)
^{4}-{r}^{3}
 \left( \sin \left( \theta \right)  \right) ^{3} \left( \cos \left(
\theta \right)  \right) ^{2}\\
& & +9\,{r}^{2} \left( \cos \left( \theta
 \right)  \right) ^{2} \left( \sin \left( \theta \right)  \right) ^{2}
+5\,{r}^{4} \left( \sin \left( \theta \right)  \right) ^{2} \left(
\cos \left( \theta \right)  \right) ^{4}-9\,r \left( \sin \left(
\theta \right)  \right) ^{3}\\
& & +28\,r \left( \cos \left( \theta \right)
 \right) ^{2}\sin \left( \theta \right) +25\,{r}^{3} \left( \cos
 \left( \theta \right)  \right) ^{4}\sin \left( \theta \right) +5\,{r}
^{5}\sin \left( \theta \right)  \left( \cos \left( \theta \right)
 \right) ^{6}\\
& & -8\, \left( \sin \left( \theta \right)  \right) ^{2}+16\,
 \left( \cos \left( \theta \right)  \right) ^{2}+20\,{r}^{2} \left(
\cos \left( \theta \right)  \right) ^{4}\\
& & +7\,{r}^{4} \left( \cos
 \left( \theta \right)  \right) ^{6}+{r}^{6} \left( \cos \left( \theta
 \right)  \right) ^{8}].
\end{array}\right.
\end{equation}

In these equations, the origin of equations (\ref{ode-uv}) is
blown up to the invariant circle $r=0$, and the complement of the
origin becomes the cylinder $r>0$. Trajectories that tend to the
origin in equation (\ref{ode-uv}) yield trajectories that tend to
an equilibrium point of equations (\ref{ode-rth}) on the circle
$r=0$. Now the zeros of $d\theta/dt$ on the circle $r=0$ are
equilibria of the rescaled equations obtained from
(\ref{ode-rth}). They are located at points where $\cos^2(\theta)
= 0,1$ or $1/3$. The equilibria determine the directions in which
trajectories of (\ref{ode-uv}) can approach or leave the origin.
These directions correspond to different approach to the
equilibrium $(0,1)$ of (\ref{ode-ac}). In the $(a,c)$ coordinates,
the directions $\theta=\pm \pi/2$ correspond to approaching
$(0,1)$ along curves tangent to the $c$-axis with tangency degree
greater than $2$. The direction $\theta=0$ corresponds to
approaching $(0,1)$ along curves tangent to $c=1$. The directions
$\pm \theta_0$ with $\cos^2\theta_0=1/3$ correspond to approaching
$(0,1)$ along curves asymptotic to the parabola $a=\frac{1}{2}
(c-1)^2$. Observe that this is consistent with Theorem
\ref{th-xcfsl2}(3). Our goal is to show that this can only happen
along a particular curve.

Since the circle $r=0$ is invariant, the Jacobians at the
equilibria discussed above are triangular. The stability of each
equilibrium is determined by the signs of $\frac{1}{r}(dr/dt)$ and
$\partial(d\theta/dt)/\partial \theta$ when these are non-zero.
The equilibria with $\cos^2(\theta) = 0$ have $\frac{1}{r}(dr/dt)
= 2$ and $\partial(d\theta/dt)/\partial \theta = -4$, so the point
is a saddle with an unstable manifold in the region $r>0$.
Equilibria with $\cos^2(\theta) = 1/3$ have $\frac{1}{r}(dr/dt) =
-4/3$ and $\partial(d\theta/dt)/\partial \theta = 16/3$, so these
equilibria are also saddles but with stable manifolds in the
region $r>0$. After change of coordinates, only one of these
stable manifolds, call it $\gamma_0$, belongs to the region
$\{c<1\}$. This curve $\gamma_0$ provides the only way to approach
$(0,1)$ which is consistent with case (3) of Theorem
\ref{th-xcfsl2}. The equilibria with $\cos^2(\theta) = 1$ have
$\frac{1}{r}(dr/dt) = 0$ and $\partial(d\theta/dt)/\partial \theta
= -8$, so further analysis is required to determine the properties
of nearby trajectories. However, because of Theorem
\ref{th-xcfsl2}, it is clear that any solution of (\ref{ode-ac})
approaching the line $\{c=1\}$ has $a\rightarrow \infty$, hence
cannot approach $(1,0)$. The following argument recovers this fact
directly from (\ref{ode-rth}).

Note that when $\theta = 0$, $d\theta/dt = 0$ and $dr/dt =
(r^5+r^7)/2$. Therefore, the $r$ axis is invariant and weakly
unstable. The trajectory along this axis approaches $(0,0)$ as $t
\to -\infty$. To prove that no other trajectories in $r>0$
approach the origin as $t \to \pm \infty$, we consider the vector
field $Y$ defined by subtracting $r^5 \cos^8(\theta)/2 +r^7
\cos^{10}(\theta)/2$ from $dr/dt$ in $X$:

\begin{equation}\label{ode-rtht}
\left\{\begin{array}{ccl} dr/dt &=& 1/2\,\sin \left( \theta
\right) r [ -5\,{r}^{3} \left( \cos
 \left( \theta \right)  \right) ^{6}-4\,r \left( \cos \left( \theta
 \right)  \right) ^{4}\\
& & -29\,{r}^{2} \left( \cos \left( \theta \right)
 \right) ^{4}\sin \left( \theta \right) -20\,\sin \left( \theta
 \right)  \left( \cos \left( \theta \right)  \right) ^{2}-35\,r
 \left( \cos \left( \theta \right)  \right) ^{2} \left( \sin \left(
\theta \right)  \right) ^{2}\\
& & -{r}^{5} \left( \cos \left( \theta
 \right)  \right) ^{8}-11\,{r}^{4} \left( \cos \left( \theta \right)
 \right) ^{6}\sin \left( \theta \right) -33\,{r}^{3} \left( \cos
 \left( \theta \right)  \right) ^{4} \left( \sin \left( \theta
 \right)  \right) ^{2}\\
& & -6\,{r}^{5} \left( \cos \left( \theta \right)
 \right) ^{6} \left( \sin \left( \theta \right)  \right) ^{2}-10\,{r}^
{4} \left( \cos \left( \theta \right)  \right) ^{4} \left( \sin
 \left( \theta \right)  \right) ^{3}\\
& & -2\,{r}^{3} \left( \cos \left( \theta \right)  \right) ^{2}
\left( \sin \left( \theta \right)
 \right) ^{4} -18\,{r}^{2} \left( \cos \left( \theta \right)  \right) ^
{2} \left( \sin \left( \theta \right)  \right) ^{3}+6\,r \left(
\sin
 \left( \theta \right)  \right) ^{4}\\
& & +2\,{r}^{2} \left( \sin \left( \theta \right)  \right)
^{5}+4\, \left( \sin \left( \theta \right)
 \right) ^{3} ],
\\
d\theta/dt &=& -1/2\,\cos \left( \theta \right) \sin \left( \theta
\right)  [ -2 \,{r}^{2} \left( \sin \left( \theta \right)  \right)
^{4}-{r}^{3}
 \left( \sin \left( \theta \right)  \right) ^{3} \left( \cos \left(
\theta \right)  \right) ^{2}\\
& & +9\,{r}^{2} \left( \cos \left( \theta
 \right)  \right) ^{2} \left( \sin \left( \theta \right)
\right) ^{2}+5\,{r}^{4} \left( \sin \left( \theta \right)  \right)
^{2} \left( \cos \left( \theta \right)  \right) ^{4}-9\,r \left(
\sin \left(
\theta \right)  \right) ^{3}\\
& & +28\,r \left( \cos \left( \theta \right)
 \right) ^{2}\sin \left( \theta \right) +25\,{r}^{3} \left( \cos
 \left( \theta \right)  \right) ^{4}\sin \left( \theta \right)\\
& &  +5\,{r}^{5}\sin \left( \theta \right)  \left( \cos \left(
\theta \right)
 \right) ^{6}-8\, \left( \sin \left( \theta \right)  \right) ^{2}
+16\,\left( \cos \left( \theta \right)  \right) ^{2}+20\,{r}^{2}
\left(
\cos \left( \theta \right)  \right) ^{4}\\
& & +7\,{r}^{4} \left( \cos
 \left( \theta \right)  \right) ^{6}+{r}^{6} \left( \cos \left( \theta
 \right)  \right) ^{8} ].

\end{array}\right.
\end{equation}
The vector field $Y$ is transverse to the vector field $X$ in the
interior of the first quadrant: $d\theta/dt < 0$ for both $X$ and
$Y$ and the $r$ component of $X$ is larger than the $r$ component
of $Y$. Therefore, trajectories of $X$ cross the trajectories of
$Y$ from below to above as they move left in the $(\theta,r)$
plane. The vector field $Y$ has a common factor of $\sin(\theta)$
in its two equations. When $Y$ is rescaled by dividing by this
factor, the result is a vector field that does not vanish in a
neighborhood of the origin. Since the $\theta$ axis is invariant
for $Y$, the $Y$ trajectory $\gamma$ starting at $(\theta,r), \,
r>0$ approaches the $r$ axis at a point with $r> 0$. The $X$
trajectory starting at $(\theta,r)$ lies above $\gamma$, so it
does not approach the origin. This proves that the only
trajectories of equations (\ref{ode-rth}) asymptotic to the origin
lie on the $r$ and $\theta$ axes.

\begin{figure} [!hbp]
   \begin{center}
   \begin{tabular}{c}
   \epsfig{figure=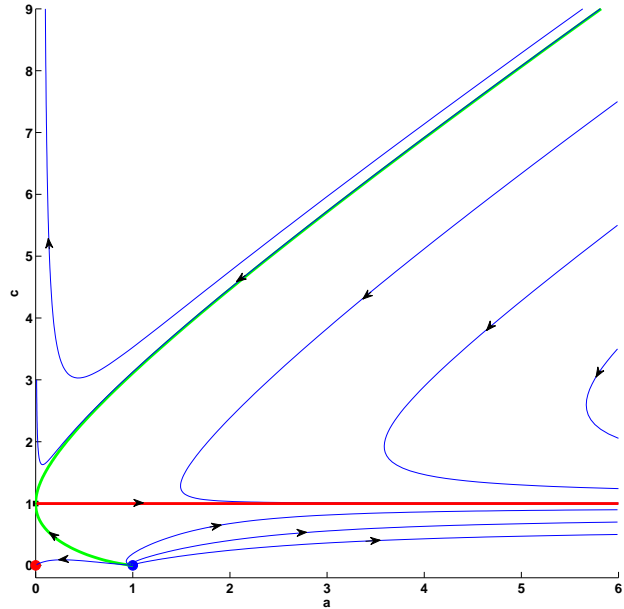,width=10cm}
   \end{tabular}
   \end{center}
   \caption
   { \label{fig:xcf}
   The flow line diagram of system (\ref{ode-ac})}
\end{figure}

In conclusion, the above analysis shows that the regions $Q_1$,
$Q_2$ are separated by the $2$-dimensional cone $S_0$ determined
by the curve $\gamma_0$ in the $(a,c)$-plane. This proves Theorem
\ref{thm3.3}. Figure \ref{fig:xcf} describes the flow lines of
(\ref{ode-ac}). The part of interest to us is the part below the
line $\{c=1\}$, which corresponds to $\{B>C\}$. The part above the
line $\{c=1\}$ corresponds to the case $\{B<C\}$, where the role
of $B$ and $C$ are exchanged. The most important component of this
diagram are the flow lines that are forward asymptotic to $(0,1)$.
They correspond to the hypersurface $S_0$ in Theorem
(\ref{thm3.3}). The flow lines in the upper-left corner have
vertical asymptotes $\{a=a_{\infty}\}$ with all positive values of
$a_{\infty}$ appearing. Similarly, the flow lines
 on the right have horizontal asymptotes $\{c=c_{\infty}\}$ with
all positive values of $c_{\infty}$ appearing. These facts can
easily be derived from the system (\ref{ode-ac}). This proves the
part of Remark \ref{rmk3.1} dealing with case (1) of Theorem
\ref{thm3.3}. The other case is similar using different
coordinates.

\bibliographystyle{halpha}
\bibliography{bio}
\end{document}